\documentclass[12pt]{article}
\usepackage{amsmath,amssymb,amsthm,xspace}

\numberwithin{equation}{subsection}

\newtheorem{thm}{Theorem}[section]

\newtheorem{prop}[thm]{Proposition}

\newcommand{\extraline}{\vskip\baselineskip}

\renewcommand{\th}{^\text{th}}

\newcommand{\R}{\mathbb{R}}
\newcommand{\Z}{\mathbb{Z}}
\newcommand{\rationals}{\mathbb{Q}}
\newcommand{\Zd}{{\Z^d}}
\newcommand{\vect}{\vec}
\newcommand{\hyp}{{H}}
\newcommand{\Tn}{T(\vect 0, \vect n)}
\newcommand{\THn}{T(\vect 0, \hyp_n)}
\newcommand{\coex}{\text{\rm Coex}}
\newcommand{\geo}{\text{\rm Geo}}
\newcommand{\sides}{\text{\rm Sides}}
\newcommand{\E}{\operatorname{\mathbf{E}}}
\renewcommand{\Pr}{\operatorname{\bf Pr}}
\newcommand{\Var}{\operatorname{\mathbf{Var}}}
\newcommand{\Q}{{\mathbf{Q}}}
\newcommand{\ind}[1]{{\mathbf{1}_{#1}}}

\newcommand{\Haggstrom}{H\"aggstr\"om\xspace}

\newcommand{\Ale}{MR1202516}
\newcommand{\BDJ}{MR1682248}
\newcommand{\BRa}{math.PR/0602496}
\newcommand{\BRb}{math.PR/0609730}
\newcommand{\BSclt}{math.PR/0412369}
\newcommand{\BKS}{MR2016607}
\newcommand{\Boi}{MR1074741}
\newcommand{\CD}{MR624685}
\newcommand{\DHa}{MR2216472}
\newcommand{\DHb}{MR2217820}
\newcommand{\DL}{MR606981}

\newcommand{\Ede}{MR0136460}
\newcommand{\GM}{MR2115045}
\newcommand{\GMdensity}{math.PR/0608667}
\newcommand{\HM}{MR1379157}
\newcommand{\Ham}{MR0405665}
\newcommand{\Hofa}{MR2114988}
\newcommand{\Hofb}{math.PR/0508114}

\newcommand{\HNeuclidean}{MR1452554}
\newcommand{\HNanimals}{MR1703127}
\newcommand{\HNgeodesics}{MR1849171}
\newcommand{\HPa}{MR1659548}
\newcommand{\HPb}{MR1794536}
\newcommand{\HW}{MR0198576}
\newcommand{\Joh}{MR1757595}
\newcommand{\JohShape}{MR1737991}
\newcommand{\JohToeplitz}{MR1957518}
\newcommand{\Kal}{MR1876169}
\newcommand{\Kesa}{MR876084}
\newcommand{\Kesb}{MR1221154}

\newcommand{\Kina}{MR0254907}
\newcommand{\Kinb}{MR0356192}
\newcommand{\KPZ}{PhysRevLett.56.889}
\newcommand{\LN}{MR1387641}
\newcommand{\LNP}{MR1421992}
\newcommand{\Lig}{MR806224}
\newcommand{\Mar}{math.PR/0301055}
\newcommand{\New}{MR1404001}
\newcommand{\NP}{MR1349159}
\newcommand{\PimMultitype}{math.PR/0411583}
\newcommand{\PimAsymptotic}{math.PR/0510605}
\newcommand{\PP}{MR1283187}
\newcommand{\Ric}{MR0329079}
\newcommand{\Ros}{MR635270}
\newcommand{\Sep}{MR1640344}
\newcommand{\Tala}{MR1303654}
\newcommand{\Talb}{MR1361756}
\newcommand{\TW}{MR1257246}
\newcommand{\Ulam}{MR0129165}
\newcommand{\VAW}{MR1094141}
\newcommand{\VAWshape}{MR1166620}

\newcommand{\WA}{MR1069633}
\newcommand{\Weh}{MR1453745}
\newcommand{\WW}{MR1617053}
\newcommand{\WR}{MR0478390}

\title{\Large First Passage Percolation and Competition Models}
\author{Nathaniel D. Blair-Stahn}
\date{April 20, 2007}

\begin{document}

\maketitle

\addcontentsline{toc}{section}{Abstract}

\begin{abstract}
This paper is a survey of various results and techniques in first
passage percolation, a random process modeling a spreading fluid on
an infinite graph. The latter half of the paper focuses on the
connection between first passage percolation and a certain class of
stochastic growth and competition models.
\end{abstract}

\newpage

\addcontentsline{toc}{section}{Contents}
\tableofcontents

\newpage

\section{Introduction}
First passage percolation is a random process on a (typically infinite) graph. Hammersley and Welsh \cite{\HW} introduced first passage percolation
as a model of fluid flow through a randomly porous medium. In this
model, each edge $e$ in the graph is assigned a random nonnegative
number $\tau(e)$, called the passage time of $e$, which is
interpreted as the time it takes to cross the edge in either
direction. In other contexts, $\tau(e)$ may represent a weight or a
capacity, but we shall stick with the passage time interpretation.
The picture to keep in mind is that of a fluid emanating from some
source vertex and flowing outward through the edges of the graph
according to the prescribed passage times. Equivalently, one can
think of an infection spreading out from some initial locus and
transmitted between neighboring vertices at random times.

First passage percolation can be defined on any graph, but the most
commonly studied model is the one in which the underlying graph is
the integer lattice $\Zd$ for $d\ge 2$. This is the model we will
focus on, although we briefly discuss models on certain random
infinite graphs in \S~\ref{isotropic sec}. The most basic results in
first passage percolation rely on ergodic theory and the theory of
subadditive processes. In fact, the study of first passage
percolation was an impetus for the development of Kingman's
subadditive ergodic theorem \cite{\Kina}. A good general reference
detailing the fundamental results in first passage percolation is
\cite{\Kesa}.

Based on heuristic arguments, the growing interface described by
first passage percolation is believed to belong to the Kardar-Parisi-Zhang (KPZ)
universality class \cite{\KPZ} of models in statistical physics. In
particular, the Eden growth model \cite{\Ede}, which can be
described in terms of a particular first passage percolation
process, has been studied extensively in this context.  We will
define the Eden growth model in \S~\ref{1-type sec}. In
\S~\ref{deviations sec} we will see some of the progress that has been made
in verifying various predictions from statistical physics.

The rest of the paper is organized as follows.
In \S~\ref{fpp sec}, after giving the precise definition of the first passage percolation
process and introducing some of the topics of interest, we summarize
the early results in the subject and proceed to describe some of the
more recent work that has been carried out. In \S~\ref{Richardson
sec} we describe how to define growth processes and competition
models based on first passage percolation and summarize recent work
in this direction.



%







\section{First passage percolation}\label{fpp sec}
\subsection{Basic definitions}
Let $\Zd$ be the integer lattice of dimension $d\ge 2$, which we consider both as a graph and as a subset of $\R^d$. Two vertices $u,v\in \Zd$ are
adjacent if $||u-v||_1 = 1$, and we denote the edge set of $\Zd$ by $E(\Zd)$. Let $\{\tau(e)\}_{e\in E(\Zd)}$ be a collection of nonnegative random variables indexed by the edges. We call $\tau(e)$ the \emph{passage time} of the edge $e$, and it represents the time needed to cross the edge in either direction. The joint law of the passage times $\tau(e)$ determines the first passage percolation process.

In order to do anything useful with the first passage percolation model, we need to make some assumptions about the distribution of the passage times.
Typically, the minimal assumption one makes is that the
passage times $\{\tau(e)\}_{e\in E(\Zd)}$ are stationary and ergodic
with respect to translations of $\Zd$. More explicitly, we can
consider the canonical sample space $\Omega=(\R_+)^{E(\Zd)}$,
equipped with some probability measure $\nu$ defined on the product
$\sigma$-field. If $\omega\in \Omega$ is a realization of
$\nu$, the passage times for $\omega$ are given by $\tau_\omega(e) = \omega(e)$.
Now, for each $u\in\Zd$, let $\theta_u : \Omega\to\Omega$ be the
natural shift operator defined by
$$
\theta_u \omega(e) = \omega(e+u),
$$
where the notation ``$e+u$" has the obvious meaning. The passage
times $\tau(e)$ are \emph{stationary} if for each $u$, the measure
$\nu$ is $\theta_u$-invariant (i.e.\ $\nu\circ\theta_u^{-1} = \nu$).
Stationary passage times are \emph{ergodic} if any event which is
invariant under every $\theta_u$ has probability 0 or 1 (i.e.\ any event $A\subset\Omega$ such that $\theta_u^{-1}A = A$ for all $u$ must satisfy $\nu(A)\in \{0,1\}$).
Throughout the rest of the paper we will assume that the passage times are stationary and ergodic, and we will be most interested in the
case where they are in fact independent and identically distributed (i.i.d.) and have finite
expectation.  In particular, our focus in \S~\ref{Richardson sec}
will be on i.i.d.\ exponential passage times. \extraline

We now define several concepts that will be discussed in more detail in later sections. Suppose that $\{\tau(e)\}_{e\in E(\Zd)}$ is a collection of passage
times describing a first-passage percolation process on $\Zd$. If
$\gamma$ is a path in $\Zd$, then the \textbf{passage time of
$\gamma$} is
$$
T(\gamma) = \sum_{e\in\gamma} \tau(e).
$$
If $U,V\subset \Zd$, the \textbf{passage time from $U$ to $V$} is
\begin{equation} \label{passage time def eqn}
T(U,V) = \inf\{T(\gamma) : \gamma \text{ is a path from $U$ to
$V$}\}.
\end{equation}
(If $U$ or $V$ is a singleton, we will write its unique element in
place of the set when using this notation or other similar notation.) We can extend this
definition to subsets of $\R^d$ as follows:  If $A\subset \R^d$, let
$\widetilde{A}$ consist of all the lattice points that are closest
to $A$, i.e.
$$
\textstyle \widetilde{A}
    = \left\{v\in \Zd : v\in x+\left[-\frac{1}{2}, \frac{1}{2}\right]^d
    \text{ for some } x\in A\right\},
$$
and for $U,V\subset \R^d$, set $T(U,V) := T(\widetilde{U},
\widetilde{V})$.

For example, for each $n=0,1,2,\ldots$, let
$$
    \vect n = (n,0,\ldots,0)\in \Zd
    \quad \text{and} \quad
    \hyp_n = \{z \in \Zd : z_1= n\}.
$$
We refer to $T(\vect 0,\vect n)$ as a \textbf{point-to-point passage
time} and $T(\vect 0,\hyp_n)$ as a \textbf{point-to-hyperplane} (or
\textbf{point-to-line} when $d=2$) \textbf{passage time}. We will
see in \S~\ref{subadditive sec} that the passage times $\Tn$ and
$\THn$ satisfy a law of large numbers (Theorem~\ref{time const
thm}), which shows that first passage percolation has an asymptotic
speed along the coordinate axes.

One of the primary objects of interest in first passage percolation
is the set $B(t)$ of vertices that can be reached from the origin by
time $t$, or a continuum version $\bar{B}(t)$ of this set in which
each $v\in B(t)$ is replaced with a unit cube centered at $v$. That
is,
$$
    B(t) = \{v\in \Zd : T(\vect 0, v) \le t\}
    \quad \text{and} \quad
    \bar B(t) = \{x\in \R^d : T(\vect 0, x) \le t\}.
$$
If we think of the percolation process as modeling an infection
spreading outward from the origin, $B(t)$ is the set of vertices
which are infected at time $t$. (We will return more explicitly to
this interpretation of $B(t)$ in \S~\ref{Richardson sec}, where we
discuss growth and competition models.) One of the fundamental
results about first passage percolation is that, under some mild
hypotheses for the passage times, $\bar B(t)/t$ converges almost
surely to a deterministic shape (Theorem~\ref{shape thm}), so the
process in fact has an asymptotic speed in all directions
simultaneously.  We will discuss this so called shape theorem
further in \S~\ref{shape sec}, and in \S~\ref{deviations sec} we
will discuss the related question of deviation bounds for the
convergence.

Another topic that arises naturally in the study of first passage
percolation is that of time-minimizing paths, or geodesics. For sets
$U,V\subset \R^d$, if $T(U,V)=T(\gamma)$ for some (necessarily
finite) lattice path $\gamma$ (i.e.\ $\gamma$ achieves the infimum
in \eqref{passage time def eqn}), we call $\gamma$ a
\textbf{geodesic from $U$ to $V$}, and we denote any such path by
$G(U,V)$. More generally, a finite or infinite path $\gamma$ in $\Zd$ is
called a \textbf{geodesic} if every finite subpath $\gamma'$ of
$\gamma$ satisfies $T(\gamma') = T(u',v')$, where $u'$ and $v'$ are
the endpoints of $\gamma'$.  Observe that a finite path $\gamma$
with endpoints $u$ and $v$ is a geodesic if and only if
$\gamma=G(u,v)$. When the passage times are i.i.d., it is easy to
see that $G(u,v)$ exists and is unique a.s.\ for each $u,v\in\Zd$ if
and only if $\tau(e)$ is a continuous random variable (see e.g.\
\cite[Lemma~8]{\WW}).  See \cite[\S~4]{\GM} for conditions
guaranteeing the existence and uniqueness of finite geodesics with
stationary passage times. We will discuss the existence of infinite
geodesics in \S~\ref{geodesics sec}.

\subsection{The subadditive ergodic theorem and the time
constant}\label{subadditive sec}

One property of the point-to-point passage times $\{T(\vect m, \vect
n)\}_{0\le m < n}$ that is immediate from the definition
\eqref{passage time def eqn} is
$$
    \Tn \le T(\vect 0, \vect m) + T(\vect m, \vect n)
    \text{ for all } 0<m<n.
$$
This motivates the following definition: A doubly indexed process
$\{X_{m,n}\}_{0\le m< n}$ is called \emph{subadditive} if $X_{0,n}
\le X_{0,m}+X_{m,n}$ for all $0<m<n$.

The main result about subadditive processes is the subadditive
ergodic theorem, which was developed by Kingman \cite{\Kina} to
study point-to-point passage times and is now a standard tool in
first passage percolation and other applications. The following
version, due to Liggett \cite{\Lig}, is an improvement on Kingman's
original result.
(Instead of (b) and (c) below, Kingman
assumed that the distribution of $\displaystyle
\left\{X_{m+k,n+k}\right\}_{0\le m<n}$ does not depend on $k$, in
which case it follows from (a) that $X_{\ell,n} \le
X_{\ell,m}+X_{m,n}$ for all $\ell<m<n$.)

\begin{thm}[Subadditive ergodic theorem \cite{\Kina}, \cite{\Lig}]\label{sub erg thm}
Suppose $X_{m,n}$, $0\le m<n$, is a family of random variables
satisfying
\begin{enumerate}
\renewcommand{\labelenumi}{{\rm (\alph{enumi})}}
\item $X_{0,n} \le X_{0,m}+X_{m,n}$ for all $0<m<n$.
\label{subadditive}

\item For each $k\ge 1$, the sequence $\displaystyle \left\{X_{nk,(n+1)k}\right\}_{n\ge 0}$
is stationary. \label{stationary seq}

\item The distribution of the sequence $\displaystyle
\left\{X_{m,m+k}\right\}_{k\ge 1}$ does not depend on $m$.
\label{same distr}

\item $\E X_{0,1}^+ < \infty$. \label{finite pos mean}
\end{enumerate}
Then
\begin{enumerate}
\renewcommand{\labelenumi}{{\rm (\roman{enumi})}}
\item $\lim_{n\rightarrow\infty}\E X_{0,n}/n = \inf_n \E X_{0,n}/n =
\gamma$ for some $\gamma \ge -\infty$. \label{means converge}

\item The limit $X = \lim_{n\rightarrow\infty} X_{0,n}/n$ exists and
is less than $+\infty$ a.s. \label{limit exists}

\item If there is some $c<\infty$ such that $\E X_{0,n}^- \ge -cn$
for all $n$, then the convergence in {\rm (ii)} also holds
in $L^1$, so $EX = \gamma$. \label{finite neg mean}

\item If the stationary sequences in {\rm (b)} are
ergodic, then $X=\gamma$ a.s. \label{ergodic}
\end{enumerate}
\end{thm}

If the random variables in Theorem~\ref{sub erg thm} are all
degenerate, then the theorem reduces to a well-known result about
subadditive functions (see e.g.\ \cite[p.\ 191]{\Kal} or
\cite{\Kinb}). On the other hand, if $\displaystyle
\left\{\xi_k\right\}_{k\ge 1}$ is a stationary sequence with $\E
|\xi_k|<\infty$, then $X_{m,n} := \xi_{m+1} + \ldots + \xi_n$
satisfies (a)--(d) and (iii) (with equality in (a), so the process
$\{X_{m,n}\}$ is in fact \emph{additive}), and Theorem~\ref{sub erg
thm} reduces to Birkhoff's ergodic theorem in this case.

If the passage times $\tau(e)$ are stationary and ergodic with
finite expectation, the random variables $X_{m,n} := T(\vect m,
\vect n)$, $0\le m < n$, are easily seen to satisfy
(a)--(d), (iii), and
(iv) of Theorem~\ref{sub erg thm}, so $\Tn/n$ converges
a.s.\ to some constant $\mu_1<\infty$. The constant $\mu_1$ is known
as the \textbf{time constant} in the direction $\vect 1$, and its
reciprocal is the the asymptotic speed of the first passage
percolation process along the coordinate axes. When the passage
times are i.i.d., it turns out that the scaled point-to-hyperplane
passage times $\THn/n$ converge to the same constant $\mu_1$. This
was first proved by Wierman and Reh \cite{\WR}, and can be deduced
from the shape theorem (Theorem~\ref{shape thm} below -- see
\cite[pp.~166-167]{\Kesa}). (Note that the process $T(\vect m,
\hyp_n)$ is not subadditive, so we cannot apply
Theorem~\ref{sub erg thm} directly.) We restate these two results
for i.i.d.\ passage times in the following theorem.

\begin{thm}[Time constant \cite{\Kina}, \cite{\WR}]\label{time const thm}
Suppose the passage times $\{\tau(e)\}_{e\in \Zd}$ are i.i.d.\ with
finite expectation. Then there is a constant $\mu_1 <\infty$ such
that
$$
    \lim_{n\rightarrow\infty} \frac{T(\vect 0,\vect n)}{n}
        = \lim_{n\rightarrow\infty} \frac{T(\vect 0,\hyp_n)}{n}
        = \mu_1 \quad \text{a.s. and in $L^1$}.
$$
\end{thm}

Finally, we mention two basic results about the time constant
$\mu_1$ in the case of i.i.d.\ passage times. First, it is easy to
see that in general $\mu_1 < \E[\tau(e)]$ (see
\cite[Theorem~4.1.9]{\HW}). Also, observe that $\mu_1=0$ corresponds
to infinite percolation speed, so that the process has superlinear
growth. There is a simple criterion for deciding when this occurs.

\begin{prop}\label{zero time const prop}
For i.i.d.\ passage times $\tau(e)$, the time constant $\mu_1$ is
nonzero if and only if $\Pr[\tau(e)=0] < p_c(\Zd)$, where $p_c(\Zd)$
is the critical value for Bernoulli bond percolation on $\Zd$.
\end{prop}

A proof of Proposition~\ref{zero time const prop} can be found in
\cite[\S~6]{\Kesa}. A heuristic argument goes as follows.  If
$\Pr[\tau(e)=0] > p_c(\Zd)$, then there is a.s.\ an infinite cluster
in $\Zd$ on which the travel time between any two vertices is zero.
It will a.s.\ take only finite time to reach this cluster from the
origin, at which point the process can head off in any direction
with infinite speed.  On the other hand, if $\Pr[\tau(e)=0] <
p_c(\Zd)$, then a.s.\ all the clusters on which infinite speed can
occur have finite size. Thus, the process can only travel a finite
distance before it has to step off one of these clusters and
accumulate some positive travel time before reaching the next
cluster. It is not too hard to show that this accumulated travel
time must with high probability increase linearly with the distance
traveled, so that the asymptotic speed is finite a.s.  The situation
at the critical value $p_c(\Zd)$ is a bit more delicate, but
Proposition~\ref{zero time const prop} shows that the asymptotic
speed in this case is infinite.

\subsection{The shape theorem}\label{shape sec}
It is natural to generalize the idea of the time constant and
consider the speed of percolation in arbitrary directions rather
than just along the coordinate axis.  In particular, for any $x\in\R^d$ with
rational coordinates, we can apply Theorem~\ref{sub erg thm} to see
that there is some constant $\mu(x)$ such that $T(\vect 0,
nx)/n\rightarrow \mu(x)$ a.s.  With this notation we have $\mu_1 =
\mu(\vect 1)$. For i.i.d.\ passage times, it is not difficult to
show that the function $\mu: \rationals^d\to [0,\infty)$ is
Lipschitz continuous and hence can be extended to all of $\R^d$, and
that the resulting function $\mu: \R^d\to [0,\infty)$ is either
identically zero or defines a norm on $\R^d$. We will refer to $\mu$
as the norm for the first passage percolation process when
appropriate; more generally, we will refer to $\mu$ as the
\textbf{shape function} for the process because of its role in the
shape theorem, which we now describe.

Recall the definitions of the growing shapes
$$
    B(t) = \{v\in \Zd : T(\vect 0, v) \le t\}
    \quad \text{and} \quad
    \bar B(t) = \{x\in \R^d : T(\vect 0, x) \le t\}.
$$
Under some moment conditions on the passage times, it can be shown
that $\bar B(t)/t$ converges almost surely to the unit $\mu$-ball
$B_0 = \{x\in \R^d : \mu(x)\le 1\}$ as $t\to\infty$. This result is
known as the shape theorem, and an in probability version was first
proved by Richardson \cite{\Ric} for $d=2$. Cox and Durrett
\cite{\CD} used a result of Kesten (found in \cite[p.~903]{\Kinb}) to strengthen Richardson's result to an almost sure
version. The following version, valid in any dimension, is proved by Kesten in
\cite{\Kesa}.
\begin{thm}[Shape theorem {\cite[Thm.~1.7]{\Kesa}}]\label{shape thm}
Suppose that $\{\tau(e)\}_{e\in E(\Zd)}$ are i.i.d.\ passage times
such that $\E \min\{\tau(e_1)^d,\ldots,\tau(e_{2d})^d\} <\infty$
(where $\{e_1,\ldots,e_{2d}\}$ is any set of $2d$ distinct edges).
Let $\mu$ be the shape function for the process, let
$\mu_1=\mu(\vect 1)$ be the time constant, and let $B_0 = \{x\in
\R^d : \mu(x)\le 1\}$.
\begin{enumerate}
\item
If $\mu_1 >0$, then $B_0$ is compact and convex with nonempty
interior, and for any $\epsilon >0$,
$$
    (1-\epsilon)B_0 \subset \frac{\bar B(t)}{t} \subset (1+\epsilon)B_0
$$
for all large $t$ almost surely.

\item If $\mu_1=0$, then $\mu\equiv 0$ (so $B_0=\R^d$), and for any
compact set $K\subset \R^d$,
$$
K \subset \frac{\bar B(t)}{t}
$$
for all large $t$ almost surely.
\end{enumerate}
\end{thm}

The moment condition on the passage times in Theorem~\ref{shape thm}
is optimal, in the sense that if it fails then $\limsup_{v\to\infty}
T(\vect 0, v)/||v||_1 =\infty$ a.s. However, even without any moment
conditions on the passage times $\tau(e)$, it is possible to define
modified passage times $\widehat{T}(u,v)$ for $u,v\in\Zd$ and a
corresponding set $\widehat{B}(t)$ such that an analogue of
Theorem~\ref{shape thm} holds (see \cite{\Kesa}). By
Proposition~\ref{zero time const prop}, we see that $B_0=\R^d$ if
and only if $\Pr[\tau(e)=0] < p_c(\Zd)$. The convexity of $B_0$
follows from subadditivity, and when $B_0\ne \R^d$, compactness and
nonempty interior follow from the fact that $\mu$ is a norm.
Otherwise, little is known about the limit shape $B_0$ other than
the obvious fact that it must have all the symmetries of $\Zd$.
Kesten \cite[\S~8]{\Kesa} shows that if the passage times are
i.i.d.\ exponential and $d$ is large, then $B_0$ is \emph{not} a
Euclidean ball, casting doubt on the conjecture that $B_0$ might be
a disc for $d=2$ based on early Monte Carlo simulations \cite{\Ede}.
Durrett and Liggett \cite{\DL} show that there are i.i.d.\ passage
times for which $B_0$ has flat edges but is not a diamond or a
square. In particular, this occurs if $\tau(e)$ is nontrivial but
attains some nonzero minimum value with probability greater than
$p_c^{\text{dir}}(\Zd)$, where $p_c^{\text{dir}}(\Zd)$ is the
critical value for \emph{directed} Bernoulli bond percolation on
$\Zd$.

There is also a version of the shape theorem for stationary passage
times. Boivin \cite{\Boi} proves that if the passage times $\tau(e)$
are stationary, ergodic, and have finite moment of order
$d+\epsilon$ for some $\epsilon>0$, then $B(t)/t$ converges a.s.\ to
a deterministic shape $B_0$. In the
stationary case, the shape function $\mu$ may take on both zero and
strictly positive values so that the limit shape $B_0$ can be an
unbounded proper subset of $\R^d$. However, if $\mu(x)>0$ for every
unit vector $x$, then $B_0$ is compact, convex, has nonempty
interior, and is symmetric with respect to reflection through the
origin. ($B_0$ may fail to have further symmetries since isotropy
may not hold in the non-i.i.d.\ case.) Conversely, \Haggstrom and
Meester \cite{\HM} show that any set $B_0\subset \R^d$ with these
properties can arise as the limit shape for some collection of
stationary passage times.

\subsection{Deviations in the passage times and the growing shape}\label{deviations sec}
Throughout this section we will assume that the passage times
$\{\tau(e)\}_{e\in E(\Zd)}$ are i.i.d.\ and satisfy the hypotheses
of Theorem~\ref{shape thm} so that $B(t)/t\to B_0$ a.s. We further
assume that $\Pr[\tau(e)=0] < p_c(\Zd)$ so that $\mu_1
>0$ and the limit shape $B_0$ is compact. \extraline

\subsubsection{The variance of $\Tn$}

Once we know that $B(t)$ converges, we can ask how much it deviates
from the limit shape $B_0$. There are various ways to approach this
problem. As a first step, we consider the variance of
$\Tn$. It is predicted that the standard deviation of $\Tn$ is of
order $n^{\chi}$ for some constant $\chi=\chi(d)$. Based on
heuristic arguments from statistical physics, it is expected that
$\chi(2)={1/3}$ (see e.g.\ \cite{KrugSpohn}, \cite{\KPZ}). This
conjecture is supported by simulations
and by rigorous results for related growth models (e.g.\
\cite{\BDJ}, \cite{\JohShape}, \cite{\Joh}), which we shall discuss in
\S~\ref{LPP sec}. The situation is less clear for higher dimensions
$d$, although it is generally believed that $\chi$ is nonincreasing
in $d$ (see \cite{\NP} for a discussion). So far, the only general
bound on $\chi$, due to Kesten, is $\chi(d)\le 1/2$ for all $d$:

\begin{thm}[Kesten \cite{\Kesb}]\label{var thm}
If $\E[\tau(e)^2] <\infty$, then there are positive constants $c_1$
and $c_2$ such that
$$
c_1  \le \Var[\Tn] \le c_2 n.
$$
\end{thm}
Kesten proves Theorem~\ref{var thm} using martingale methods (the
``method of bounded differences"). Although Theorem~\ref{var thm}
provides the best known bounds for a general distribution on the
passage times, better bounds have been proved for certain classes of
distributions.  For example, Benjamini, Kalai, and Schramm
\cite{\BKS} use an inequality of Talagrand \cite[Thm.~1.5]{\Tala}
to show that $\Var [\Tn] = O(n/\log n)$ if the passage
times have the uniform distribution on $\{a,b\}$, where $0<a<b$. They note
that the essential feature of first passage percolation needed to
prove both their result and Kesten's is that the number of edges
$e\in E(\Zd)$ such that modifying $\tau(e)$ increases $\Tn$ is
bounded by a constant times $n$. Building on the methods in
\cite{\BKS}, Bena\"im and Rossignol \cite{\BRa} use a Gaussian
version of Talagrand's \cite{\Tala} inequality and apply the
techniques of \cite{\BKS} to prove $O(n/\log n)$ variance for a
large class of i.i.d.\ absolutely continuous passage times,
including exponential.

As for lower bounds on the variance,  Newman
and Piza \cite{\NP} prove that in dimension $d=2$,
$\Var[\Tn]=\Omega(\log n)$ under certain hypotheses on the passage
times. In particular, if we set $\lambda = \inf \{x : \Pr(\tau(e)
\le x) >0\}$, the condition assumed in \cite{\NP} is that
\begin{equation}\label{nonflat eqn}
\Pr(\tau(e)=\lambda) < p(\lambda),
\end{equation}
where
$$
p(\lambda) =
\begin{cases}
    p_c(\Z^2) & \text{if } \lambda=0,\\
    \vspace{5pt}
    p_c^{\text{dir}}(\Z^2) & \text{if } \lambda>0.
\end{cases}
$$
Based on Proposition~\ref{zero time const prop} and the results in
\cite{\DL}, this condition is necessary for the shape $B_0$ to be
compact and for its boundary to have no flat edges; it is suspected
that \eqref{nonflat eqn} should also be sufficient for this to hold
(see \cite{\NP}). Pemantle and Peres \cite{\PP} use different
methods to prove $\Omega(\log n)$ variance for the special case of
exponential passage times in $d=2$. On the other hand, for higher
dimensions $d$, it is still not known whether the variance of $\Tn$ even
diverges as $n\to\infty$. \extraline

\subsubsection{Large deviation bounds for $\Tn$}

If the assumption of finite variance for the passage times is
strengthened to the existence of a finite exponential moment, then
one can obtain good bounds on the deviation of $\Tn$ from its
expected value, and on the deviations of $\E[\Tn]$ from $n\mu$.  The
following theorem is due primarily to Kesten \cite{\Kesb}, with the
upper bound in (\ref{ETn - n mu eqn}) being an improvement made by
Alexander.
\begin{thm}[Kesten \cite{\Kesb}, Alexander \cite{\Ale}]\label{deviations
thm}

If $\E [e^{\gamma \tau(e)}] < \infty$ for some $\gamma >0$, then
there exist positive constants $c_1,c_2,c_3,c_4,c_5$, such that
\begin{equation}\label{tailbounds eqn}
    \Pr\left(\left|\frac{\Tn -\E[\Tn]}{\sqrt n}\right|\ge x\right)
    \le c_1 e^{-c_2 x} \quad \text{for } x\le c_3 n,
\end{equation}
and
\begin{equation}\label{ETn - n mu eqn}
c_4\frac{1}{n} \le \E[\Tn] - n\mu \le c_5 n^{1/2} \log n.
\end{equation}
\end{thm}
Note that the lower bound in \eqref{ETn - n mu eqn} strengthens the
trivial inequality $\E[\Tn] \ge n\mu$ implied by Theorem~\ref{sub
erg thm}. Both Theorems~\ref{var thm} and Theorem~\ref{deviations
thm} remain valid if $\Tn$ is replaced by $\THn$ (see \cite{\Kesb}
or \cite{\Ale}), or if $\vect n$ is replaced by any $v\in\Zd$ and
$n$ is replaced by $||v||_1$. In fact, using versions of
\eqref{tailbounds eqn} and \eqref{ETn - n mu eqn} valid for
arbitrary directions, Kesten \cite[Theorem~2]{\Kesb} shows that there is some
constant $C$ (depending on the dimension $d$ and the distribution of $\tau(e)$) such that
almost surely,
\[
\left(1 -  \left(\frac{C\log t}{\sqrt t}\right)^{\frac{1}{d+2}}\right)\cdot B_0
\subset \frac{B(t)}{t} \subset
\left(1 + \frac{C\log t}{\sqrt t}\right)\cdot B_0
\quad\text{for all large } t.
\]
Some improvements of Theorem~\ref{deviations thm} are available in
certain situations. Talagrand \cite[\S 8.3]{\Talb} shows that the
upper bound in \eqref{tailbounds eqn} can be strengthened to
$O(e^{-cx^2})$ if $\E[\Tn]$ is replaced with a median of $\Tn$.
For the same class of distributions considered in \cite{\BRa} (with
the added assumption of finite exponential moment), Bena\"im and
Rossignol \cite{\BRb} prove that \eqref{tailbounds eqn} still holds
if the $\sqrt n$ in the denominator is replaced by $\sqrt{n/\log
n}$. Instead of the Talagrand-type inequalities used in \cite{\BRa},
the techniques used in \cite{\BRb} involve modified Poincar\'e
inequalities arising from the context of ``threshold phenomena" for
Boolean functions. \extraline

\subsubsection{Scaling exponents for the growth process}

In the statistical physics literature (see e.g.\ \cite{KrugSpohn}), the
fluctuations of a randomly growing shape such as $\bar B(t)$ are
studied in terms of two exponents $\chi$ and $\xi$, which describe
respectively the longitudinal and transverse fluctuations in the
surface of $\bar B(t)$. For example, it is expected that the
standard deviation of the time $T(\vect 0, H)$ at which $\bar B(t)$
first reaches a hyperplane $H$ at distance $r$ from the origin is of
order $r^\chi$, while the set of points in $H$ which are likely to
be first reached by $\bar B(t)$ is expected to have diameter on the
order of $r^\xi$. There are various ways to define $\chi$ and $\xi$
precisely, and it is an open problem to determine whether the various
definitions are equivalent.

The exponents $\chi$ and $\xi$ are not
expected to depend on the underlying distribution of the
$\tau(e)$'s, at least under certain hypotheses (for example,
\eqref{nonflat eqn} above -- see \cite{\NP} or \cite{\LNP}). A
priori, $\chi$ and $\xi$ could depend on the direction of
travel, but it is expected that they should be the same in any
direction in which the boundary of $B_0$ has nonzero curvature, at least in low dimensions. The values of
$\chi$ and $\xi$ \emph{are} expected to depend on the dimension $d$,
but heuristic arguments suggest that the scaling identity $\chi=2\xi-1$ holds in
all dimensions (see \cite{KrugSpohn}). As noted in the introduction, first passage percolation models are expected to belong to the KPZ universality class \cite{\KPZ}, leading to the prediction that $\chi(2)=1/3$ and (in accordance with the scaling identity) $\xi(2)=2/3$.

We now describe some of the progress that has been made towards
computing the exponents $\chi$ and $\xi$. Since $\chi$ and $\xi$
might depend on the direction of travel, we will write $\chi_{\hat
x}$ and $\xi_{\hat x}$ to denote their values in the direction of
some unit vector $\hat x\in\R^d$. In \cite{\NP}, Newman and Piza
show that in any dimension $d$, if \eqref{nonflat eqn} holds, then
$\chi_{\hat x} \ge (1-(d-1)\xi_{\hat x})/2$ (this was proved by Wehr
and Aizenman \cite{\WA} for $d=2$). Then they show that, under the
same hypothesis \eqref{nonflat eqn}, if the passage times have
finite exponential moment, then $\xi_{\hat x}\le 3/4$ for any $\hat
x$ which is a direction of curvature for $B_0$ (i.e.\ a direction in
which the boundary of $B_0$ has nonzero curvature). For $d=2$, this
yields $\chi_{\hat x}\ge 1/8$ in any direction of curvature $\hat
x$, improving the previously mentioned logarithmic lower bound on
$\Var[T(\vect 0, n\hat x)]$. It is easy to show that any compact
convex set has a direction of curvature \cite[Lemma~5]{\NP}, so in
$d=2$ there is at least one direction $\hat x$ such that
$\Var[T(\vect 0, n\hat x)] = \Omega(n^{1/4})$ when the $\tau(e)$'s
have finite exponential moment.

The method in \cite{\NP} used to prove $\xi_{\hat x}\le 3/4$ makes use of an exponent $\chi'$ analogous
to $\chi$, but which also takes into account the deviations of
$\E[T(\vect 0, n\hat x)]$ from $n\mu(\hat x)$. The Kesten-Alexander
deviation bounds (Theorem~\ref{deviations thm}) imply that $\chi'\le
1/2$. Newman and Piza then use a rigorized version of the heuristic
argument from \cite{KrugSpohn} for the scaling identity $\chi=2\xi-1$ to show
that $\xi_{\hat x} \le (1+\chi')/2$, which yields the bound
$\xi_{\hat x}\le 3/4$.

In \cite{\LNP}, Licea, Newman, and Piza extend the methods in
\cite{\NP} to obtain lower bounds on various versions of the
exponent $\xi$. Combining the trivial bound $\chi \ge 0$ with the
(nonrigorous) scaling identity yields the nontrivial bound $\xi\ge
1/2$, which is expected to hold in all dimensions. The value
$\xi=1/2$ corresponds to what is called a diffusive process, and it
is believed that, at least in low dimensions, first passage
percolation should in fact be superdiffusive, i.e.\ $\xi>1/2$. Using
progressively weaker definitions $\xi^{(1)}$, $\xi^{(2)}$,
$\xi^{(3)}$ for $\xi$, Licea, Newman, and Piza prove
$$
    \xi^{(1)}(d)\ge1/(d+1), \quad
    \xi^{(2)}(d)\ge 1/2, \quad
    \text{and}\quad \xi^{(3)}(2)\ge 3/5,
$$
assuming that the passage times satisfy \eqref{nonflat eqn} and/or
$\E[\tau(e)^2]<\infty$. The latter two bounds correspond to superdiffusivity as predicted by the physical models. While the first bound is subdiffusive, it is nontrivial from a mathematical perspective, and may be useful because the exponent $\xi^{(1)}$ has certain advantages over the other two definitions of $\xi$.

\subsection{Infinite geodesics} \label{geodesics sec}
Recall that a geodesic is a time-minimizing path in first passage percolation, and that $G(U,V)$ denotes a geodesic between the sets $U$ and $V$ when such a path exists.
Suppose that $G(u,v)$ exists and is unique for each pair of
vertices $u,v\in\Zd$. For any $u\in\Zd$, we define the \textbf{tree
of infection of $u$}, $\Gamma(u)$, to be the (graph theoretic) union
of all the finite geodesics starting at $u$:
$$
    \Gamma(u) = \bigcup_{v\in\Zd}G(u,v).
$$
The fact that $\Gamma(u)$ is a tree follows from the uniqueness of
the geodesics. If we think of the percolation process as modeling an
infection spreading outward from $u$, then the unique path in
$\Gamma(u)$ from $u$ to another vertex $v$ traces the route by which
$v$ became infected.

Let $K(\Gamma(u))$ denote the number of topological ends in
$\Gamma(u)$ -- that is, the number of semi-infinite paths in
$\Gamma(u)$ starting at $u$. We call any such path a
\textbf{one-sided geodesic} starting at $u$. A standard compactness
argument shows that $K(\Gamma(u))\ge 1$ for any $u$. (The set of finite geodesics starting at $u$ can be viewed in a natural way as a compact space, so it must contain a limit point since it has infinitely many elements.) In \cite{\New},
Newman uses the Kesten-Alexander deviation bounds (Theorem~\ref{deviations thm}) and methods similar to those in \cite{\NP} to show that if the passage times are i.i.d.\ and the curvature
of the boundary of $B_0$ is uniformly bounded away from 0, then
$K(\Gamma(u)) =\infty$ a.s.\ for any $u$. While the assumption of
uniform curvature is plausible, there are no i.i.d.\ probability
measures on the passage times for which $B_0$ is known to have this
property. Hoffman \cite{\Hofb} has shown that $K(\Gamma(u))=\infty$
under a much weaker assumption on the limit shape, namely that $B_0$
is not a polygon. Although there are no i.i.d.\ passage times which
are known to satisfy this assumption either, the result of
\cite{\HM} shows that there are stationary passage times for which
it holds.  We will revisit this topic in \S~\ref{equal rates sec},
where we discuss the connection between the existence of geodesics
and the question of coexistence in a certain competition model
obtained as a projection of first passage percolation.

Newman proves the above result by considering geodesics with an asymptotic direction. If $\hat x$ is a unit vector in $\R^d$ and $\gamma$ is a one-sided geodesic with vertices $v_0,v_1,v_2,\ldots$, then $\gamma$ has asymptotic direction $\hat x$ if $\lim_{n\to\infty}v_n/||v_n||_2 = \hat x$, and we call $\gamma$ an  \textbf{$\hat x$-geodesic}. It is not known in general whether $\hat x$-geodesics exist or whether every one-sided geodesic must have a direction, but Newman \cite{\New} gives affirmative answers to both questions under the assumption that $B_0$ is uniformly curved. Licea and Newman \cite{\LN} use some of the ideas in \cite{\New} to prove a uniqueness result for $\hat x$-geodesics when $d=2$:

\begin{thm}[Licea and Newman \cite{\LN}]\label{coalesce thm}
Suppose the passage times $\{\tau(e)\}_{e\in\Z^2}$ are i.i.d.\ with continuous distribution. Then for Lebesgue almost every $\hat x$ on the unit circle,
$$
    \Pr\left(\text{\rm There exist disjoint $\hat x$-geodesics}\right)=0,
$$
and hence any two $\hat x$-geodesics must coalesce.
\end{thm}

Naturally, in addition to one-sided geodesics, one can also consider
\textbf{two-sided geodesics}, i.e.\ geodesics which are infinite in
both directions instead of just one. In contrast to one-sided geodesics, it is unclear whether two-sided geodesics exist at all. (Since a two-sided geodesic does not have a fixed starting point, the compactness argument used to prove the existence of one-sided geodesics fails in this case.) Most of the work on two
sided-geodesics has focused on the case of continuous i.i.d.\ passage times in $d=2$, where it is expected
that two-sided geodesics do not exist. We summarize some results in
this direction.

In \cite{\Weh}, Wehr shows that almost surely, the number of 2-sided
geodesics is either 0 or $\infty$ in $d=2$, and that an analogous
result holds in $d$ dimensions for locally weight-minimizing
hypersurfaces instead of curves. This result is equivalent to the
statement that the number of ground states in the random exchange
Ising model (REIM) is 2 or $\infty$ a.s.

Using Theorem~\ref{coalesce thm}, Licea and Newman \cite{\LN} show that for Lebesge-a.e.\ unit vector $\hat x\in \R^2$, there cannot exist an $(\hat x, -\hat x)$-geodesic, i.e.\ a two-sided geodesic with asymptotic directions $\hat x$ and $-\hat x$.  Wehr and Woo \cite{\WW} show that if $H$ is any half-plane in $\R^2$, then there can exist no two-sided geodesics contained entirely within $H$. As a corollary, any two-sided geodesic must intersect every line with rational slope.

\subsection{Isotropic models of first passage percolation}\label{isotropic sec}

As we have seen on various occasions above, one disadvantage of the
first passage percolation model on $\Zd$ is that we do not have much
information about the limit shape $B_0$. One way to get around this
is to define a stochastically isotropic model so that symmetry
considerations imply that $B_0$ must be a Euclidean ball.  For
example, Vahidi-Asl and Wierman \cite{\VAW} introduce models in
which the underlying graph is either a random Vorono\u\i\
tesselation of the plane or its dual Delaunay triangulation, where
the centers of the Vorono\u\i\ cells are given by a Poisson point
process on $\R^2$. Howard and Newman \cite{\HNeuclidean} introduce a
different model, in which the underlying graph is the complete graph
with vertices given by a Poisson point process on $\R^d$ and the
passage times are given by $\tau(e)=|e|^\alpha$, where $\alpha>1$
and $|e|$ denotes the Euclidean distance between the endpoints of
$e$.

Using a random graph for the for the percolation process introduces
various technical problems, but nevertheless, versions of many of
the results familiar from the $\Zd$ model still hold for these
Euclidean models. For example, in the Vorono\u\i\ and Delaunay
models, there is a time constant \cite{\VAW}, a shape theorem
\cite{\VAWshape}, and deviation bounds similar to those in
Theorems~\ref{var thm} and \ref{deviations thm}
\cite{\PimAsymptotic}. Furthermore, since $B_0$ has uniform
curvature in this model, Pimentel \cite{\PimMultitype} is able to
use the techniques in \cite{\NP} and \cite{\New} to show that the
transversal fluctuation exponent $\xi\le 3/4$, and that almost
surely, every one-sided geodesic has an asymptotic direction and
there exists a one-sided geodesic in every direction. Similar
results are proved by Howard and Newman for their model in
\cite{\HNeuclidean}, \cite{\HNanimals}, and \cite{\HNgeodesics}. We
also mention that \cite{\PimMultitype} contains results about a
competion model on the Delaunay triangulation analogous to the
competition model on $\Zd$ described in \S~\ref{equal rates sec}
below.


\subsection{Directed first passage and last passage percolation}
\label{LPP sec}

The process we have been referring to as first passage percolation
is more properly called \emph{undirected} first passage percolation.
One can also consider \emph{directed} first passage percolation or a
related model called (directed) last passage percolation. Both of
these models are defined similarly to undirected first passage
percolation, except that only \emph{increasing} paths (defined
below) are allowed. Certain versions of directed last passage
percolation are much better understood than undirected first passage
percolation. We now describe the directed models and summarize some
of the most interesting results.

\extraline

Let $\{\tau(e)\}_{e\in\Zd}$ be a collection of i.i.d.\ nonnegative
passage times.  A path $\gamma$ in $\Zd$ is called
\textbf{increasing} if each step in $\gamma$ is made by increasing a
single coordinate by 1. For $u,v\in\Zd$, write $u\le v$ if $u_i\le
v_i$ for $1\le i\le d$. If $u\le v$, define the (directed)
\textbf{first-passage time from $u$ to $v$} to be
$$
    T_{\min}(u,v) = \min \{T(\gamma) : \gamma
    \text{ is an increasing path from $u$ to $v$}\},
$$
and define the \textbf{last-passage time from $u$ to $v$} to be
$$
    T_{\max}(u,v) = \max \{T(\gamma) : \gamma
    \text{ is an increasing path from $u$ to $v$}\},
$$
where $T(\gamma)$ is defined as in the undirected case. As before,
we can extend these definitions to passage times between two points
in $\R^d$. The directed models are often defined with passage times
$\tau(v)$ on the vertices $v$ of $\Zd$ rather than the edges, but
the analysis is similar with either convention, so we will stick
with edge passage times.

As with the undirected first-passage times, the directed
first-passage times are subadditive, whereas the last-passage times
are \emph{superadditive}, i.e.\ for vertices $u\le v\le w$ we have
$$
    T_{\max}(u,v) + T_{\max}(v,w) \le T_{\max}(u,w).
$$
Applying Theorem~\ref{sub erg thm} to the first-passage times and a
superadditive version of Theorem~\ref{sub erg thm} to the last-passage times implies that there
are shape functions $g,h:(\R_+)^d\to[0,\infty)$ such that, for all
$x\in(\R_+)^d$,
$$
    \lim_{n\rightarrow\infty} \frac{T_{\min}(\vect 0, nx)}{n}= g(x)
    \text{\ \ a.s.\ \quad and \quad }
    \lim_{n\rightarrow\infty} \frac{T_{\max}(\vect 0, nx)}{n}= h(x)
    \text{\ \ a.s.}
$$
Furthermore, we can define growing shapes analogous to $B(t)$:
$$
    U(t) = \{x\in (\R_+)^d : T_{\min}(\vect 0, x)\le t\}
    \quad \text{and}\quad
    V(t) = \{x\in (\R_+)^d : T_{\max}(\vect 0, x)\le t\}.
$$
Under appropriate conditions on the distribution of $\tau(e)$,
Martin \cite{\Mar} proves a shape theorem for the directed models:
$$
    {U(t)}/{t} \to U_0 \text{\ \ a.s.\ \quad and \quad}
    {V(t)}/{t} \to V_0 \text{\ \ a.s.},
$$
where $U_0=\{x : g(x)\le 1\}$ and $V_0=\{x : h(x)\le 1\}$. In the first-passage case, subadditivity implies that $U_0$ is convex, whereas in the last-passage case, superadditivity implies that $(\R_+)^d\setminus V_0$ is convex.

In contrast with the undirected model, there are two special cases
of directed last-passage percolation in $d=2$ for which the shape
function $h(x)$ is known explicitly. A theorem of Rost \cite{\Ros}
(see also \cite{\BSclt}, \cite{\Mar}) implies that for exponential
passage times with mean 1,
$$
    h(x) = ||x||_{1/2} = (\sqrt{x_1} + \sqrt{x_2}\,)^2.
$$
Johansson \cite{\JohShape} shows that for geometric passage times
with parameter $q$,
\begin{equation*}\label{geo shape eqn}
    h(x) = h_q(x) = \frac{q(x_1+x_2) + 2\sqrt{x_1 x_2}}{1-q}.
\end{equation*}
These are the only two nontrivial cases where the shape function for
i.i.d.\ passage times is known, in any of the directed or
undirected, first- or last-passage models. However, Sepp\"al\"ainen
\cite{\Sep} finds the limiting shape for a particularly simple
stationary model of directed first-passage percolation on
$(\Z_+)^2$, in which vertical edges have a deterministic, constant
passage time and horizontal edges have i.i.d.\ Bernoulli passage
times.

In fact, Johansson \cite{\JohShape} not only identifies the shape function in the i.i.d.\ geometric last-passage model, but extends the techniques of Baik, Deift, and Johansson \cite{\BDJ}
to show that the passage times $T_{\max}(\vect 0, n x)$, appropriately centered and scaled, converge in distribution to the Tracy-Widom \cite{\TW} distribution for the largest eigenvalue in a random matrix sampled from the Gaussian Unitary Ensemble (GUE).
In particular, it is shown that the standard deviation of $T_{\max}(\vect 0, n x)$ in this model is of order $n^{1/3}$, so that $\chi=1/3$ in accordance with the predictions of KPZ universality \cite{\KPZ}.

\extraline

Finally, we mention another model which can be viewed as a continuum version of the directed last-passage percolation model defined above, and in fact can be obtained as a limit of last-passage percolation with i.i.d.\ geometric passage times (see, e.g.\ \cite{\JohToeplitz}). This model was introduced by Hammersley \cite{\Ham} as a method for approaching Ulam's problem \cite{\Ulam} of finding the distribution of the longest increasing subsequence in a random permutation.

Consider a unit-rate Poisson process on $\R^2$. Analogous to
increasing lattice paths, we can define an increasing path between
Poisson points to be a path $\gamma$ that moves only up and to the
right.  That is, $\gamma$ is a sequence of Poisson points such that
if $x$ and $x'$ are consecutive points in $\gamma$, then $x\le x'$.
We define the length of $\gamma$ to be the number of Poisson points
it contains. Let $L(r)$ be the length of the longest increasing path
between Poisson points contained in the square $[0,r]^2$.
Conditional on the event that $[0,r]^2$ contains $N$ points, $L(r)$
has the same distribution as the length of the longest increasing
subsequence in a random (uniform distribution) permutation of
$\{1,2,\ldots,N\}$ (see \cite{\Ham}). Baik, Deift, and Johansson
\cite{\BDJ} show that
$$
\lim_{r\rightarrow\infty}\Pr\left(\frac{L(r) -2r}{r^{1/3}}\le s\right) =F(s),
$$
where F(s) is the Tracy-Widom distribution for the largest eigenvalue of a GUE random matrix. The $r^{1/3}$ in the denominator shows that $\chi=1/3$ for this model, where $\chi$ is the exponent  describing the longitudinal fluctuations of a maximal increasing path. Moreover, Johansson \cite{\Joh} applies the techniques from \cite{\NP} and \cite{\LNP} to show that the transversal fluctuations of the maximal paths have exponent $\xi=2/3$, verifying the scaling identity $\chi=2\xi-1$ for this model.

\section{Richardson's growth model and competition models}
\label{Richardson sec}

\subsection{The 1-type Richardson model}\label{1-type sec}
First passage percolation can be put into the framework of interacting particle systems (see e.g. \cite{MR776231}) by defining a $\{0,1\}^\Zd$-valued process $\{\eta_t\}_{t\ge 0}$ given by
\begin{align*}
\eta_t(v) &=
\begin{cases}
    1 & \text{if } v\in B(t), \\
    0 & \text{otherwise},
\end{cases}
&  \text{for } v\in\Zd.
\end{align*}
We may think of sites in state 0 as healthy and sites in state 1 as infected, so that the process represents an infection spreading outward from the origin.

When the passage times are i.i.d.\ exponentials with parameter
$\lambda>0$, the memoryless property implies that the process $\eta_t$ is Markovian. In this case the process is called
\textbf{Richardson's growth model} \cite{\Ric}, also known as the ``contact
process with no recoveries," in comparison with the similarly defined contact process (cf.\ \cite{MR776231}).  In Richardson's growth model, a site in state 1
remains infected forever and tries to infect each of its $2d$
neighbors at rate $\lambda$. Thus, the rate at which an uninfected
site flips from 0 to 1 is equal to $\lambda$ times the number of
infected neighbors it has. Only the origin is infected at time 0.

As noted in \cite{\Ric}, this process is related to a discrete time
process called Eden's growth model \cite{\Ede}, defined as follows:
Set $A_1=\{\vect 0\}$, and for $n>1$ set $A_n = A_{n-1} \cup
\{v_n\}$, where $v_n$ is chosen from the set of uninfected sites
with probability proportional to the number of neighbors it has in
$A_{n-1}$. Then $A_n$ has the same distribution as $B({t_n})$,
where
$$
t_n = \inf \{t :  B(t) \text{ contains $n$ vertices}\}.
$$

\subsection{The 2-type Richardson model}\label{2-type sec}
In \cite{\HPa}, H\"aggstr\"om and Pemantle introduced the
\textbf{two-type Richardson model}. In this model, instead of one
type of particle spreading throughout the lattice, there are two
species of particles competing for space. This competition is
described by a $\{0,1,2\}^\Zd$-valued Markov process
$\{\xi_t\}_{t\ge 0}$ with parameters $\lambda_1$ and $\lambda_2$
which determine the flip rates as follows: 1's and 2's never flip,
while a 0 flips to a 1 (resp. a 2) at rate $\lambda_1$ (resp.
$\lambda_2$) times the number of neighbors of type 1 (resp. 2). The
1's and 2's represent sites infected by species 1 and 2,
respectively, and 0's represent uninfected sites.

One natural question we can ask in the two-type model is whether
both species continue growing indefinitely or whether one species
ends up surrounded by the other so that it is only able to infect a
finite number of sites.  If $A_1$ and $A_2$ are two disjoint subsets
of $\Zd$, we denote by $\coex(A_1, A_2)$ the event that both species
eventually infect an infinite number of sites when species $i$
initially occupies the sites in $A_i$, and we call this event
\textbf{coexistence} or \textbf{mutual unbounded growth} for the
initial configuration $(A_1, A_2)$. It is easy to see that
$\Pr(\coex(A_1, A_2))<1$ unless both of the sets $A_1$ and $A_2$ are
already infinite, so the first nontrivial question to ask is whether $\Pr(\coex(A_1,
A_2)) >0$. Clearly coexistence is impossible if one of the sets
$A_i$ surrounds the other set $A_j$, i.e.\ if there is no infinite
path starting in $A_j$ that does not intersect $A_i$. We say that
the pair $(A_1, A_2)$ is \textbf{fertile} if neither set surrounds
the other. Deijfen and H\"aggstr\"om showed that as long as the
initial configuration of the process is finite and fertile, the
choice of configuration is irrelevant to the question of whether
coexistence has positive probability:
\begin{thm}[Deijfen and H\"aggstr\"om \cite{\DHa}] \label{init config thm}
If $(A_1, A_2)$ and $(A_1', A_2')$ are two fertile pairs of disjoint
finite sets in $\Zd$, then for any pair of growth rates $\lambda_1$
and $\lambda_2$,
$$
\Pr(\coex(A_1, A_2))>0 \Leftrightarrow \Pr(\coex(A_1', A_2'))>0.
$$
\end{thm}
Weaker versions of Theorem~\ref{init config thm} (in which the sets $A_i$ and $A_i'$ consist of single points) appeared in \cite{\HPa} and \cite{\GM}, and most treatments of coexistence have simply focused on the case where the initial configuration is $(\{\vect 0\},\{\vect 1\})$.

Intuitively, if the growth rates $\lambda_1$ and $\lambda_2$ are equal, we might expect that coexistence occurs with positive probability since neither species has an inherent advantage over the other. This result was in fact proved for $d=2$ in \cite{\HPa}, and subsequently generalized to any $d\ge 2$ in \cite{\GM} and \cite{\Hofa}. On the other hand, if the growth rates are different, say $\lambda_1> \lambda_2$, then unless species 2 gets lucky and surrounds species 1 relatively quickly, species 1 is likely to overtake species 2 by virtue of its superior speed, making coexistence implausible.  H\"aggstr\"om and Pemantle conjecture in \cite{\HPa} that $\Pr(\coex(\vect 0,\vect 1)) = 0$ when $\lambda_1\ne \lambda_2$ and prove a somewhat weakened version of this conjecture in \cite{\HPb}. In the next two sections we discuss the two-type Richardson model in more detail in the two cases $\lambda_1= \lambda_2$ and $\lambda_1\ne \lambda_2$.

\subsection{Competition with equal growth rates}\label{equal rates sec}

The two-type Richardson model is somewhat simpler to analyze when
both species grow at the same rate. If $\lambda_1= \lambda_2=
\lambda$, we can obtain the two-type process $\xi_t$ as the
projection of a single first passage percolation process with
i.i.d.\ exponential($\lambda$) passage times, analogous to the
definition of $\eta_t$ in the one-type model. If the process starts
with initial configuration $(A_1,A_2)$, then
$$
\xi_t(v) =
\begin{cases}
    1 & \text{if } T(A_1,v) \le t \text{ and } T(A_1,v) < T(A_2,v)\\
    2 & \text{if } T(A_2,v) \le t \text{ and } T(A_2,v) < T(A_1,v)\\
    0 & \text{otherwise}.
\end{cases}
$$
Note that the definition of $\xi_t$ can be generalized in an obvious
way to model competition between $k$ species with equal growth rates
and initial configuration $(A_1,\ldots,A_k)$, for any $k\ge 1$.
Since the two-type (or $k$-type) model and the one-type model are
both defined in terms of an underlying first passage percolation
process, results about one model can often be translated into
results about the other, as will be illustrated below.

H\"aggstr\"om and Pemantle first addressed the question of
coexistence for species with equal growth rates in \cite{\HPa},
where they proved that $\Pr[\coex(\vect 0, \vect 1)]>0$ when $d=2$.
The main step in their proof was to show that in the related one-type process starting at $\vect 0$, there are
infinitely many sites in the right half-plane which have a $>50\%$
probability of being infected after their neighbor to the left, so that these sites ``sense" that the infection is coming from the left. From there, it is a small step to show that in the two-type process, with positive probability there are infinietly many sites in the right half plane that are reached by species 2 before they are reached by species 1, and that a symmetric situation holds in the left half-plane.

Observe that the definition of $\xi_t$ makes sense for more general
passage time distributions, although the Markov property holds only
in the i.i.d.\ exponential case. However, for any stationary
distribution of passage times, the two species will still be growing
at the same average rate, and we might expect coexistence to hold in
the stationary case as well.
Indeed, Garet and Marchand \cite{\GM} and Hoffman \cite{\Hofa}
independently generalized the coexistence result of \cite{\HPa} to a
large class of stationary ergodic passage times in any dimension
$d\ge 2$. Furthermore, an analogue of Theorem~\ref{init config thm}
holds in the stationary case so that the starting configuration is
still irrelevant \cite[p.~312]{\GM}.

Coexistence in the two-type or $k$-type model is related to the
existence of one-sided geodesics in the corresponding one-type
model. If coexistence of $k$ species occurs, then the same
compactness argument used to show that $K(\Gamma(\vect 0))\ge 1$
shows that there exist $k$ disjoint one-sided geodesics in the
underlying first passage percolation process, one starting in each
of the initial sets $A_1,\ldots, A_k$. Therefore, denoting
coexistence in the $k$-type model by $\coex(A_1,\ldots, A_k)$ and
the existence of disjoint geodesics $G_i$ starting in the sets $A_i$
by $\geo(A_1,\ldots, A_k)$, we have
\begin{equation}\label{geo>=coex eqn}
\Pr[\geo(A_1,\ldots, A_k)] \ge \Pr[\coex(A_1,\ldots, A_k)].
\end{equation}
Furthermore, if $\geo(A_1,\ldots, A_k)$ occurs, it seems plausible that some finite modification of passage times might allow the construction of $k$ one-sided geodesics starting at $\vect 0$ so that $\Pr[K(\Gamma(\vect 0))\ge k] >0$.  In fact, at least when $k=2$ (and probably for any $k$ -- see \cite{\Hofb}), we can go in the other direction as well, from geodesics to coexistence: For the class of stationary measures considered in
\cite{\GM} or \cite{\Hofa}, it can be shown (see
\cite[Lemma~5.3]{\GM}) that
$$
    \Pr[\coex(\vect 0, \vect 1)]>0 \Leftrightarrow
    \Pr[K(\Gamma(\vect 0))\ge 2] >0.
$$
Thus, since coexistence of two species occurs with positive
probability, there are at least two one-sided geodesics starting at
$\vect 0$ with positive probability. In fact, while Garet and
Marchand \cite{\GM} use techniques similar to those in \cite{\HPa}
to prove that coexistence is possible and then conclude that there
are at least two geodesics with positive probability, Hoffman
\cite{\Hofa} first proves that there almost surely exist at least
two distinct one-sided geodesics (not necessarily with the same
starting point) and uses this to show that coexistence has positive
probability.

In \cite{\Hofb}, Hoffman applies the techniques in \cite{\Hofa} to
the $k$-type model to obtain further results about both geodesics
and coexistence when $d=2$ and the passage times are given by a
certain class of ``good" ergodic stationary measures $\nu$. Although
the results are stated only for $d=2$, the methods can be applied to
any $d\ge 2$. We now state the main results, which depend on the
geometry of the limit shape $B_0$ corresponding to $\nu$. For a good
measure $\nu$, let $\sides(\nu)$ be the number of sides of $\partial
B_0$ if $\partial B_0$ is a polygon or infinity if $\partial B_0$ is
not a polygon.
\begin{thm}[Hoffman \cite{\Hofb}]\label{k coex thm}
Let $\nu$ be a good stationary measure on $(\R_+)^{E(\Z^2)}$, and
let $k\le\sides(\nu)$. For any $\epsilon>0$, if $r$ is sufficiently
large there exist $u_1,\ldots,u_k\in \widetilde{\partial (rB_0)}$
such that
$$
    \Pr[\coex(u_1,\ldots, u_k)] > 1-\epsilon
    \text{\ and (by \eqref{geo>=coex eqn}) }
    \Pr[\geo(u_1,\ldots, u_k)] > 1-\epsilon.
$$
\end{thm}
The points $u_i$ in Theorem~\ref{k coex thm} are chosen to be the
lattice points closest to points $u_1',\ldots,u_k'\in \partial
(rB_0)$ at which the tangent lines of $\partial (rB_0)$ are
distinct (such points exist by the assumption that
$k\le\sides(\nu)$).  The idea of the proof is that if
$v_1,\ldots,v_k$ are points on $\partial B_0$ with distinct tangent
lines $L_{v_1},\ldots,L_{v_k}$, then with positive probability, for
each $i$ there will be infinitely many $n$ such that $v_i$ is closer
(in travel time) to the translated line $nv_i + L_{v_i}$ than any of
the other points $v_j$ are.  This shows that if the process starts
with initial configuration $(v_1,\ldots,v_k)$ (assuming
$v_i\in\Z^2$), each $v_i$ will infect infinitely many sites with
positive probability, so coexistence occurs. By scaling the picture
up by a sufficiently large factor $r$, the probability of
coexistence can be made arbitrarily close to 1.

Hoffman also obtains the following results about one-sided geodesics
starting at $\vect 0$.
\begin{thm}[Hoffman \cite{\Hofb}]\label{nu geodesics thm}
Let $\nu$ be a good stationary measure on $(\R_+)^{E(\Z^2)}$.  If
$k\le \sides(\nu)/2$, then
$$K(\Gamma(\vect 0))\ge k \text{\ \ a.s.}$$
\end{thm}

\begin{thm}[Hoffman \cite{\Hofb}]\label{exp geodesics thm}
Let $\nu=(\mathcal{L}(\tau))^{\otimes E(\Z^2)}$, where $\tau$ is an
exponential random variable. If $k\le \sides(\nu)$, then
$$\Pr[K(\Gamma(\vect 0))\ge k]>0.$$
\end{thm}

Observe that by symmetry, we must have $\sides(\nu)\ge 4$ when
$d=2$, so Theorem~\ref{k coex thm} shows that coexistence of four
species is possible for any good measure $\nu$, and Theorem~\ref{exp
geodesics thm} shows that with i.i.d.\ exponential passage times it
is possible to get four one-sided geodesics starting at $\vect 0$.
Furthermore, in \cite{\DL} it is shown that there is a nontrivial
i.i.d.\ measure $\nu$ such that $B_0$ is neither a square nor a
diamond, so by symmetry we must have $\sides(\nu)\ge 8$. Thus,
Theorem~\ref{k coex thm} implies that there is a nontrivial i.i.d.\
measure $\nu$ for which coexistence of eight species is possible. In
\cite{\HM} it is shown that there exists a good measure $\nu$ such
that $B_0$ is the unit disc, so Theorem~\ref{nu geodesics thm}
implies that there exists a good measure $\nu$ such that
$K(\Gamma(\vect 0)) = \infty$ a.s.

\subsection{Competition with different growth rates}\label{different rates sec}
If the growth rates $\lambda_1$ and $\lambda_2$ are different for
the two species, we can construct the Markov process $\xi_t$ from
two independent i.i.d.\ exponential first passage percolation
processes on $\Zd$, one with parameter $\lambda_1$ and the other
with parameter $\lambda_2$.
However, the description of $\xi_t$ is not quite as
simple as it was in the case of equal growth rates because, with two
underlying sets of passage times instead of just one, there is no
guarantee that the geodesics between infected vertices of one
species will not cross geodesics of the other species. For this
reason, the values of $\xi_t$ must be defined iteratively by
considering the process at the time $t_n$ of the $n\th$ infection.
We mention that the state of the process at time $t_n$ can be
described analogously to Eden's growth model for the one-type
process, except that the edge causing the next infection is chosen
from all edges on the boundary with probability proportional to
$\lambda_i$ if the edge borders species $i$.

In \cite{\HPa}, H\"aggstr\"om and Pemantle conjecture that
coexistence in the two-type Richardson model is impossible if
$\lambda_1\ne \lambda_2$. While the full conjecture is still an open
problem, H\"aggstr\"om and Pemantle were able to prove the slightly
weaker result that coexistence is impossible for almost all choices
of parameter values:

\begin{thm}[H\"aggstr\"om and Pemantle \cite{\HPb}]\label{noncoex
thm} For the two-type Richardson model on $\Zd$, $d\ge 2$, with
$\lambda_1=1$ we have
$$
\Pr(\coex(\vect 0,\vect 1))=0
$$
for all but at most countably many choices of $\lambda_2$.
\end{thm}

By time scaling, the probability of coexistence depends only on the
ratio $\lambda = \lambda_2/\lambda_1$, so Theorem~\ref{noncoex thm}
remains true for any other choice of $\lambda_1$.  Furthermore, the
same time-scaling argument plus symmetry implies that the
probabilities of coexistence for the pairs $(1,\lambda)$ and
$(1,1/\lambda)$ are equal, so it suffices to consider the case
$\lambda_1=1$ and $\lambda_2=\lambda\in [0,1]$.

At first glance, it may seem strange that Theorem~\ref{noncoex thm}
has not been extended to include all values of $\lambda_2\ne
\lambda_1$. Intuitively, we expect that $\Pr(\coex(\vect 0,\vect
1))$ should decrease as $\lambda=\lambda_2/\lambda_1$ moves farther
away from 1. Since Theorem~\ref{noncoex thm} implies that we can
choose $\lambda$ arbitrarily close to 1 such that $\Pr(\coex(\vect
0,\vect 1))=0$, such monotonicity would imply that coexistence is
impossible for all $\lambda\ne 1$. However, it is not obvious how to
prove that the probability of coexistence is monotone in $\lambda$.
In fact, although this monotonicity property is certainly plausible
for the integer lattice $\Zd$, Deijfen and H\"aggstr\"om \cite{\DHb}
have shown that there are other (highly non-symmetric) graphs where
monotonicity does not hold.

We now give a brief outline of the proof of Theorem~\ref{noncoex
thm} from \cite{\HPb}. The main tool is the following proposition,
which we state as it appears in $\cite{DHunbounded}$. Let
$\mathbf{P}_{A_1,A_2}^{\lambda_1, \lambda_2}$ denote the law of the
two-type process with rates $\lambda_1$, $\lambda_2$ and initial
configuration $(A_1,A_2)$. For $i=1,2$, let $G_i$ be the event that
species $i$ finally infects an infinite number of sites (so
$\coex(A_1,A_2) = G_1\cap G_2$), and let $B_0$ denote the limit
shape for the one-type Richardson model with rate 1. Then
\begin{prop}[{\cite[Prop.~2.2]{\HPb}}, {\cite[Prop.~5.2]{DHunbounded}}]
\label{nonsurvival prop} For any $\lambda<1$ and $\epsilon>0$ we
have
$$
    \lim_{r\to\infty} \sup_{A_1,A_2}
    \mathbf{P}_{A_1,A_2}^{1, \lambda}(G_2) = 0,
$$
where the supremum is over all initial configurations $(A_1,A_2)$
such that
\begin{align}
    &A_2 \text{ is contained in } rB_0 \text{, while}\nonumber \\
    &A_1 \text{ is \emph{not} contained in } (1+\epsilon)rB_0.
    \label{head start eqn}
\end{align}
\end{prop}

For example, Proposition~\ref{nonsurvival prop} says that if we
start the process with with the slow species occupying the entire
$\mu$-ball of radius $r$ and the fast species occupying a single
site outside the $\mu$-ball of radius $(1+\epsilon)r$ (where $\mu$
is the norm for the unit rate Richardson model), the survival
probability of the slow species goes to zero as $r\to\infty$. Using
Proposition~\ref{nonsurvival prop} and the strong Markov property,
\Haggstrom and Pemantle show that if $\coex(\vect 0,\vect 1)$
occurs, the set of sites infected by both species, scaled by $t$,
converges a.s.\ to the limit shape of the slow species.

To prove Theorem~\ref{noncoex thm}, \Haggstrom and Pemantle first
describe a coupling $\mathbf{Q}$ of the processes $\mathbf{P}_{\vect
0,\vect 1}^{1, \lambda}$ ($\lambda\in [0,1]$) such that $\Q$-a.s.,
for all $t$, the set 1's at time $t$ decreases with $\lambda$, and the set of 2's
at time $t$ increases with $\lambda$. Writing $\coex(\lambda)$ for
the event that $\coex(\vect 0, \vect 1)$ occurs at parameter
$\lambda$ under the law $\Q$, they use the above result to show that
$\Q$-a.s., $\coex(\lambda)$ occurs for at most one $\lambda\in
[0,1]$. That is,
$$
    \Q\left(\ind{\coex(\lambda)}=0
        \text{ for all but at most one }
        \lambda\in [0,1]\right)=1,
$$
so by Fubini's theorem,
$$
    \sum_{\lambda\in [0,1]} \Q(\coex(\lambda))
    = \E_\Q  \sum_{\lambda\in [0,1]} \ind{\coex(\lambda)}
    \le \E_\Q 1 = 1.
$$
Therefore, since the sum on the left is finite, there can be only
countably many $\lambda$'s such that $\Q(\coex(\lambda))>0$, which
proves Theorem~\ref{noncoex thm} since $\Q(\coex(\lambda)) =
\mathbf{P}_{\vect 0,\vect 1}^{1, \lambda}(\coex(\vect 0,\vect 1))$.

While the results in \cite{\HPb} apply only to finite initial
configurations, Deijfen and H\"aggstr\"om \cite{DHunbounded}
recently proved some interesting results about coexistence in the
case when one of the initial sets $A_i$ is infinite. In particular,
they considered the cases where $A_1$ is either the hyperplane
$H=\hyp_0 = \{z\in \Zd : z_1=0\}$ (minus the origin) or the half
line \mbox{$L=\{z\in \Zd : z_1\le 0 \text{ and } z_i=0 \text{ for
all } i\ne 1\}$} (minus the origin), and $A_2=\{\vect 0\}$. Their
main result is
\begin{thm}[Deijfen and H\"aggstr\"om \cite{DHunbounded}]
\label{unbounded coex thm} For the two-type Richardson model in
$d\ge 2$ dimensions,
\begin{enumerate}
\item $\Pr[\coex(H\setminus \{\vect 0\}, \vect 0)] > 0$ if and
only if $\lambda_1 < \lambda_2$. \label{hyperplane coex}

\item $\Pr[\coex(L\setminus \{\vect 0\}, \vect 0)] >0$ if and only
if $\lambda_1\le \lambda_2$. \label{halfline coex}
\end{enumerate}
\end{thm}
The fact that coexistence is impossible if $\lambda_1 > \lambda_2$
for either $H$ or $L$ follows from Proposition~\ref{nonsurvival
prop} since whenever $A_1$ is infinite while $A_2$ is finite, the
pair $(A_1,A_2)$ satisfies \eqref{head start eqn} for all
sufficiently large $r$. The ``if" direction of \ref{hyperplane coex}
is proved by combining a shape theorem for the one-type process
starting from $H$ (which is proved using the large deviation bounds
of Theorem~\ref{deviations thm}) with a shape theorem for the
``hampered" one-type process starting from $\vect 0$ and restricted
to a cylinder about the first coordinate axis (which follows from a
standard modification of the proof of the ordinary shape theorem).
The strategy of proof is to show that when $\lambda_1 < \lambda_2$,
there is a positive probability that species 2 gets a big enough
head start over species 1 that it is able to take over the entire
cylinder without interference. The corresponding result for $L$
follows from the result for $H$ because (a rotation of) $L$ is a
subset of $H$, and the probability of survival for either species is
monotone with respect to the starting configuration
\cite[Lemma~3.1]{DHunbounded}.

For the critical case $\lambda_1 = \lambda_2$, the proof that
coexistence is possible when $A_1=L\setminus \{\vect 0\}$ follows
techniques similar to those used in \cite{\HPa}. In fact, the
coexistence result for $A_1=L\setminus \{\vect 0\}$ when $\lambda_1
= \lambda_2$ allows an easy proof of the the coexistence result in
\cite{\HPa} (see \cite[Theorem~6.1]{DHunbounded}).
The proof that coexistence is impossible when $A_1=H\setminus
\{\vect 0\}$ is rather more involved, and we will not discuss it
here. We mention, however, that this result shows that in a one-type
process started from $H$, almost surely every vertex in $H$ will
infect only a finite number of vertices in $\Zd$.

As was the case with equal growth rates, the definition of the
process $\xi_t$ in terms of first passage percolation makes sense
for more general passage times, although again, Markovity will be
lost in the non-exponential case. In \cite{\GMdensity}, Garet and
Marchand extend the results of \cite{\HPb} to include i.i.d.\
passage times which are not necessarily exponential but for which
the passage time distributions for the two species are
stochastically comparable.
In this setting, they show that for any $d$, if the slow species
survives, the fast species cannot occupy a very high density of
space (for example, ``it could not be observed by a medium
resolution satellite").  For $d=2$, they show that almost surely,
one species must finally occupy a set of full density in the plane
while the other species occupies only a set of null density. They
also obtain deviation bounds similar to those in \cite{\HPb} showing
that if coexistence occurs then the infected region in the two-type
process must grow essentially according to the law of the first
passage percolation process governing the slow species. Finally,
they prove an analogue of Theorem~\ref{noncoex thm} for families of
stochastically comparable passage times indexed by a continuous
parameter, showing that coexistence cannot occur except perhaps for
a countable set of parameters.

\addcontentsline{toc}{section}{References}

\bibliographystyle{halpha}
\bibliography{fpp}

\end{document}